\documentclass[12pt, leqno]{article}

\usepackage{amsmath}
\usepackage{amsfonts}
\usepackage{amssymb}
\usepackage{amsfonts}
\usepackage{comment}
\usepackage{graphicx}

\usepackage{amsthm}

\swapnumbers
\theoremstyle{plain}
\newtheorem{thm}[equation]{Theorem}
\newtheorem{lem}[equation]{Lemma}
\newtheorem{prop}[equation]{Proposition}
\newtheorem{cor}[equation]{Corollary}

\theoremstyle{definition}

\theoremstyle{remark}
\newtheorem{rem}[equation]{Remark}

\numberwithin{equation}{section}

\def\be{\begin{equation}}
\def\ee{\end{equation}}

\begin{document}

\def\sof{\hfill\rule{2mm}{2mm}}
\def\SS{\mathcal S}
\def\qq{{\bold q}}
\def\mn{\text{-}}

\title{\textbf{Packing sets of patterns}}


\author{
\textbf{Alexander Burstein}\\
Department of Mathematics\\
Iowa State University\\
Ames, IA 50011-2064, USA\\
\texttt{burstein@math.iastate.edu} \and
\textbf{Peter H\"ast\"o}\\
Department of Mathematics\\
Norwegian University of Science and Technology\\
7491 Trondheim, Norway\\
\texttt{peter.hasto@helsinki.fi}}

\date{\today}



\maketitle


\begin{abstract}
Packing density is a permutation occurrence statistic which
describes the maximal number of permutations of a given type that
can occur in another permutation. In this article we focus on
containment of sets of permutations. Although this question has
been tangentially considered previously, this is the first
systematic study of it. We find the packing density for various
special sets of permutations and study permutation and
pattern co-occurrence.
\end{abstract}


\section{Introduction}

The string $413223$ contains two subsequences, $133$ and $122$, each
of which is \emph{order-isomorphic} (or simply \emph{isomorphic})
to the string $122$, i.e.\ ordered in the same way as $122$. In
this situation we call the string $122$ a \emph{pattern}.
Herb Wilf first proposed the systematic study of pattern
containment in his 1992 address to the SIAM meeting on Discrete
Mathematics. However, several earlier results on pattern
containment exist, for example, those by Knuth \cite{Knuth} and
Tarjan \cite{Tarjan}.

Most results on pattern containment actually deal with
\emph{pattern avoidance}, in other words, enumerate or consider
properties of strings over a totally ordered alphabet which avoid
a given pattern or set of patterns.
There is considerably less research on other aspects of pattern
containment, specifically, on packing patterns into strings over a
totally ordered alphabet, but see
\cite{AAHHS, Hasto, HSV, Price, Stromquist} for the permutation
case and \cite{Ba, BSV, BHM} for the more general pattern case.

Although several of the above cited papers have defined packing
density for sets of patterns, virtually all of them have subsequently
restricted the attention to the case when the set contains only
one pattern. In this paper we take the first systematic step
in studying the set packing question: in Section~\ref{layerSect} we study the
packing density of so-called layered permutations which have been the
focus of much research also in the single permutation case.
In Section~\ref{maxSect} we compare the packing density of
a pair of patterns to the densities of the individual patterns
in a measure which we call covariance. In Section~\ref{avgSect}
we study the same question for average occurrence of patterns,
in which case the covariance is the actual, statistical covariance.

\subsection*{Notation}

Let $[k]=\{1,2,\dots,k\}$ be our canonical totally ordered
alphabet on $k$ letters, and consider the set $[k]^n$ of
$n$-letter words over $[k]$. We say that a pattern $\pi\in[l]^m$
\emph{occurs} in $\sigma\in[k]^n$, or $\pi$ \emph{hits} $\sigma$,
or that $\sigma$ \emph{contains} the pattern $\pi$, if there is a
subsequence of $\sigma$ order-isomorphic to $\pi$.

Given a word $\sigma\in[k]^n$ and a set of patterns
$\Pi\subseteq[l]^m$, let $\nu(\Pi,\sigma)$ be the total number of
occurrences of patterns in $\Pi$ ($\Pi$-patterns, for short) in
$\sigma$. Obviously, the largest possible number of
$\Pi$-occurrences in $\sigma$ is $\binom{n}{m}$, when each
subsequence of length $m$ of $\sigma$ is an occurrence of a
$\Pi$-pattern. Define
\[
\begin{split}
\mu(\Pi,k,n)&=\max\{\,\nu(\Pi,\sigma)\,|\,\sigma\in[k]^n\},\\
d(\Pi,\sigma)&=\frac{\nu(\Pi,\sigma)}{\binom{n}{m}}\ {\rm and}\\
\delta(\Pi,k,n)&=\frac{\mu(\Pi,k,n)}{\binom{n}{m}}=
\max\{\,d(\Pi,\sigma)\,|\,\sigma\in[k]^n\},
\end{split}
\]
the maximum number of $\Pi$-patterns in a word in
$[k]^n$, the probability that a subsequence of $\sigma$ of length
$m$ is an occurrence of a $\Pi$-pattern, and the maximum such
probability over words in $[k]^n$, respectively.

We want to consider the asymptotic behavior of
$\delta(\Pi,k,n)$ as $n\to\infty$ and $k\to\infty$.
R.\ Barton \cite{Ba} recently proved that
$$ \lim_{n\to \infty}\lim_{k\to \infty} \delta(\Pi, n,k) =
\lim_{k\to \infty}\lim_{n\to \infty} \delta(\Pi, n,k), $$
and so we can mend the definition from \cite{BHM} and
define the common limit to be the packing density
of set of the patterns $\Pi$.


\section{Sets of layered permutations}\label{layerSect}

In this section we deal with sets of layered permutations. Recall that a
permutation is said to be layered if it is a strictly increasing
sequence of strictly decreasing substrings. These substrings are called the
layers of the permutation. For instance, $21543$ is layered with layers
$21$ and $543$.

It has been shown that if $\Pi$ consists of layered permutations, then
there is a $\Pi$-maximal permutation which is layered \cite[Theorem~2.2]{AAHHS}.
For the single layered permutation case Price showed that
$$ \delta(\pi) = {m \choose m_1, \ldots, m_r}
\sup \sum_{i_1< \ldots < i_r} \lambda_1^{i_1} \cdots
\lambda_r^{i_r}, $$ where $\pi$ has $r$ layers of length $m_1$,
\ldots, $m_r$, $m$ is the total length of $\pi$ and the supremum
is taken over all partitions of unity $(\lambda_i)$. If the above
supremum is achieved with a partition with exactly $r$ parts, then
we call the permutation \textit{simple} \cite{Hasto}. The next
results shows that \cite[Theorem~3.3]{Hasto} generalizes to the
case of sets of permutations in some cases.

\begin{prop}\label{simpleProp}
Let $S$ be a set of layered permutations of length $m$ and $r$
layers such that the minimizing sequence is increasing.
Let $m^-$ be the shortest
layer of any permutation in $S$. If $\log_2(r+1)\le m^-$, then
$S$ is simple, and the packing density $\delta(S)$ equals
$$ \sup \sum_{\pi\in S} {m\choose m_1^\pi, \ldots, m_r^\pi}
\lambda_1^{m_1^\pi}\cdots \lambda_r^{m_r^\pi}, $$
where $m_i^\pi$ is the $i^{\rm th}$ layer of $\pi$ and the supremum is
taken over partitions of unity $(\lambda_i)_{i=1}^r$.
\end{prop}

\begin{proof} As in the proof of Theorem~3.3, \cite{Hasto}, we conclude
that the minimizing sequence of any $\pi\in S$ has $r$ layers. Hence
$S$ is simple. The last formula follows directly from this.
\end{proof}

One might think that it is always the case that a set of simple permutations
is simple. We have not been able to prove it, however.

In some cases it is easy to show that the condition of the
previous proposition holds. A layered permutation is said to be
increasing, if its layer sizes are increasing.
If $S$ is a set of increasing layered permutations, then
the minimizing sequence $(\lambda_i)$ is also increasing.
The proof of this fact is the same as in the case of only a
single permutation, see \cite[Lemma~3.2]{Hasto}. Another obvious
case is when the set $\Pi$ is symmetric, in the sense that it contains
all the permutations with certain layer sizes, like the set
$\{[2,1,1], [1,2,1], [1,1,2]\}$.

Let us next consider some special sets of layered permutations.
The prototypical case for the next theorem are permutations
$[m,2]$ and $[m,1,1]$. In this case the permutations differ only
in that the last two letters are interchanged, but as can be seen below, this
is not the reason that we are able to calculate the packing density.

\begin{thm}\label{l1Thm} Let $m,n\ge 2$ and let $S(m,n)$ be the set of all
permutations whose first layer has length $m$ and whose subsequent layers
have total length $n$. Then we have
$$ \delta( S(m,n) ) =
{m+n-1\choose n} \frac{(m-1)^{m-1} n^n}{(m+n-1)^{m+n-1}}. $$
Note that $\delta( S(m,n) ) = \delta([m-1,n])$ if $m\ge 3$.
\end{thm}

\begin{proof} There exists a (possibly infinite) sequence
$(\lambda_i)_{i=1}^\infty$ with $\sum \lambda_i = 1$ such that
\[
\begin{split}
 \delta( S(m,n) )
& = {m+n\choose n} \sum_{i<j} \lambda_i^m \lambda_j^n +
\sum_{p=1}^{n-1} {m+n\choose m, n-p, p }
\sum_{i<j<k} \lambda_i^m \lambda_j^{n-p}\lambda_k^p
+ \ldots \\
& = {m+n\choose n} \left( \sum_{i<j} \lambda_i^m \lambda_j^n +
\sum_{p=1}^{n-1} {n \choose p}
\sum_{i<j<k} \lambda_i^m \lambda_j^{n-p}\lambda_k^p + \ldots\right) \\
& = {m+n\choose n} \sum_{i} \lambda_i^m \left(\sum_{i<j}\lambda_j\right)^n.
\end{split}
\]
Let us denote $\Lambda_i = \sum_{i<j} \lambda_j$.
We define $\lambda_1' = c\lambda_1$ and $\lambda_i' = d\lambda_i$ for $i>1$.
Moreover, we choose $d(c)=(1-c\lambda_1)/(1-\lambda_1)$, in order to have
$\sum \lambda_i' = 1$. Since the original sequence $(\lambda_i)_{i=1}^\infty$
was maximal, we have
$$ \sum_{i} \lambda_i^m \Lambda_i^n \ge
\sum_{i} \lambda_i'^m \left(\sum_{j>i}\lambda_j'\right)^n
= c^m d^n \lambda_1^m \Lambda_1^n +
d^{m+n} \sum_{i>1} \lambda_i^m \Lambda_i^n. $$
Let us denote
$$ \alpha = \lambda_1^m \Lambda_1^n = \lambda_1^m (1-\lambda_1)^n $$
and $ \beta = \sum_{i>1} \lambda_i^m \Lambda_i^n$.
Then our previous conclusion implies that the function
$$ F(c) = c^m d(c)^n \alpha + d(c)^{m+n} \beta $$
has a maximum at $c=1$. Differentiating this function and
evaluating at $c=d=1$ gives
$$ F'(1) = \Big(m-\frac{n\lambda_1}{1-\lambda_1}\Big) \alpha - (m+n)
\frac{\lambda_1}{1-\lambda_1} \beta. $$
Since $1$ is a maximum, the derivative equals zero, so
$$ \beta = \frac{m(1-\lambda_1) - n\lambda_1}{(m+n)\lambda_1} \alpha. $$
Therefore
$$ \sum_{i} \lambda_i^m \Lambda_i^n = \alpha + \beta
= \frac{m}{(m+n)\lambda_1} \alpha =
\frac{m}{m+n} \lambda_1^{m-1} (1-\lambda_1)^n. $$
Clearly, the last expression is maximized by $\lambda_1=(m-1)/(m+n-1)$.
Therefore we have
\[
\begin{split}
\delta( S(m,n) )
& =
{m+n\choose n} \frac{m}{m+n} \sup_{0<\lambda_1<1} \lambda_1^{m-1}
(1-\lambda_1)^n \\
& =  {m+n\choose n} \frac{m}{m+n} \frac{(m-1)^{m-1} n^n}{(m+n-1)^{m+n-1}},
\end{split}
\]
as claimed.
\end{proof}

\begin{rem} We can make the previous theorem slightly more
general by allowing different first terms. Let $M\ge 3$ and
let $M>m_1> \ldots> m_r \ge 2$. Then we can find the packing density of
the set
$$ S(m_1, M-m_1)\cup \ldots \cup S(m_r, M-m_r) $$
by finding the maximum over $\lambda_1$ of the real valued function
$$ \frac{1}{M} \sum_{i=1}^r m_i \lambda_1^{m_i-1} (1-\lambda_1)^{M-m_i}. $$
The proof of this fact is very similar to the proof of Theorem~\ref{l1Thm},
and is thus omitted.
\end{rem}

Using the method of the previous proof, we get an upper bound for the packing
density of much more general types of permutations. In the general case,
the upper bound is not attained, however.

\begin{cor} Let $m_1\ge 2$. Then
$$ \delta([m_1, \ldots, m_r]) \le \frac{m_1}{m}
\frac{(m_1-1)^{m_1-1} (m-m_1)^{m-m_1}}{(m-1)^{m-1}} \delta([m_2, \ldots, m_r]), $$
where $m=m_1+ \ldots + m_r$.
\end{cor}

\begin{proof}
There again exists a sequence $(\lambda_i)$ such that
$$ \delta([m_1, \ldots, m_r]) =
{m \choose m_1, \ldots, m_r} \sum_{i_1<\ldots< i_r} \lambda_{i_1}^{m_1} \cdots
\lambda_{i_r}^{m_r}. $$
We split the sum into two parts,
$$ \alpha =  \lambda_1^{m_1} \sum_{i_2<\ldots< i_r} \lambda_{i_2}^{m_2} \cdots
\lambda_{i_r}^{m_r}, $$
and the rest, denoted by $\beta$.
As in the previous proof, we set $\lambda_1' = c\lambda_1$ and
$\lambda_i' = d\lambda_i$ for $i>1$, construct the function $F(c)$,
calculate the derivative, and set it equal to zero.
As before, we calculate
$$ \alpha + \beta = \frac{m_1}{m \lambda_1} \alpha. $$
Using a rescaling and the definition of packing density we find
$$ \frac{\alpha}{\lambda_1^{m_1}} = \sum_{i_2<\ldots< i_r} \lambda_{i_2}^{m_2} \cdots
\lambda_{i_r}^{m_r} \le (1-\lambda_1)^{m-m_1} \delta([m_2, \ldots, m_r]). $$
Therefore
\[
\begin{split}
\delta([m_1, \ldots, m_r])
& =  {m \choose m_1, \ldots, m_r} \frac{m_1}{m \lambda_1} \alpha \\
& \le  \sup_{0<\lambda_1<1}
{m \choose m_1, \ldots, m_r} \frac{m_1}{m}
\lambda_1^{m_1-1} (1-\lambda_1)^{m-m_1} \delta([m_2, \ldots, m_r]).
\end{split}
\]
Clearly the last supremum is reached for $\lambda_1 = (m_1-1)/(m-1)$,
from which the claim follows.
\end{proof}


\section{Maximal pattern co-occurrence}\label{maxSect}

We have the following obvious estimates for the packing density
of a set of two patterns:
$$ \max\{\delta(\pi_1), \delta(\pi_2)\}
\le \delta(\pi_1,\pi_2) \le \delta(\pi_1) + \delta(\pi_2). $$
We want to measure how close $\delta(\pi_1,\pi_2)$ is to
these extremes, so it makes sense to consider the ratio
$$ \frac{\delta(\pi_1) + \delta(\pi_2) - \delta(\pi_1,\pi_2)}
{\delta(\pi_1) + \delta(\pi_2) - \min\{\delta(\pi_1), \delta(\pi_2)\} } =
\frac{\delta(\pi_1) + \delta(\pi_2) - \delta(\pi_1,\pi_2)}
{\max\{\delta(\pi_1), \delta(\pi_2)\} } $$
as a measure of co-occurrence of the two patterns
$\pi_1$ and $\pi_2$.
Let us denote this ratio by $cov(\pi_1, \pi_2)$.

Using Proposition~\ref{simpleProp} we get the following simple corollary:

\begin{cor} Let $b>a\ge 2$. Then
$$ \delta([a,b],[b,a]) = {a+b\choose a}
\sup_{x\in [0,1]} x^a(1-x)^b + x^b(1-x)^a. $$
If $b-a=1$ or $b-a=2$ then the supremum occurs at
$x=1/2$.
\end{cor}

\begin{proof}
The formula for the packing density follows since
we may rearrange the layers in the optimal permutation
as we want, since the set $[a,b],[b,a]$ is symmetric.
The claim about the supremum follows by direct calculation
of the derivative.
\end{proof}

Using the previous corollary
we get the following table for the co-occurrence of
the layered permutations $[a,b]$ and $[b,a]$:

\begin{center}
\begin{tabular}{lrrrr}
$a\setminus b$ & 3 & 4 & 5 & 6 \\
2 & 0.191 & 0.576 & 0.915 &  \\
3 & -- & 0.138 & 0.447 & 0.799 \\
4 & -- & -- & 0.108 & 0.365
\end{tabular}
\end{center}

Let us next calculate the packing density for sets of
patterns of length three. For permutations this was done in \cite{AAHHS}.

\begin{prop} We have
$$ \delta(112) = \delta(112, 121) = \delta(112,121,211) = 2\sqrt{3} - 3. $$
\end{prop}

\begin{proof}
The numerical value $\delta(112) = 2\sqrt{3} - 3$ is from \cite[Example~2.12]{BHM}.
We next complete the proof by showing that $\delta(112) = \delta(112,121,211)$.
The remaining equality follows from this, since the density
certainly grows if we add more permutations to a set.

Let $\sigma$ be a word and
consider adjacent distinct letters at $\sigma_i$ and $\sigma_{i+1}$
and let $\sigma'$ be the pattern with these letters interchanged.
Then
$$ d(112,121,211; \sigma) = d(112,121,211; \sigma'). $$
To see this notice that the number of occurrences of $112$ which hit
at most one of the two letters at position $i$ and $i+1$ is the
same in $\sigma$ and $\sigma'$. The same holds for the other two
patterns.
So it remains to consider occurrences involving both of these
positions. Assume $\sigma_i<\sigma_{i+1}$. Then if $112$ hits
$\sigma$ at positions $j<i<i+1$ it is clear that $121$ hits
$\sigma'$ in the same positions. Similarly a hit of $121$ at
$i<i+1<j$ is turned into a hit of $211$ at the same positions. If
$\sigma_i>\sigma_{i+1}$, then the situation is reversed. Hence in
each case the total number of occurrences is preserved.

We have now shown that we may exchange adjacent letters in
$\sigma$. Doing this an sufficient number of times we may assume
that $\sigma$ is increasing. But then all the hits are of type $112$,
hence
$$ d(112,121,211; \sigma) = d(112; \sigma). $$
Since $\sigma$ was arbitrary, the claim follows.
\end{proof}

\begin{prop}
We have
$\displaystyle \delta(112, 122) = \delta([2,1],[1,2]) = \frac{3}{4}$.
\end{prop}

\begin{proof}
Since both $112$ and $122$ are non-decreasing,
it is clear that the minimizing pattern must be non-decreasing. We
may assume that the minimizing pattern is of the form
$$ \sigma = 1^{s_1}\, 2^{s_2}\, \ldots\, n^{s_n}. $$
Consider then the permutation of type $\sigma'=[s_1, \ldots,
s_n]$. It is clear that every occurrence of $112$ in $\sigma$
corresponds to an occurrence of $[2,1]$ in $\sigma'$, similarly
for $122$ and $[1,2]$.
\end{proof}

\begin{cor}
We have $\displaystyle
\delta(112,121,211, 221, 212, 122) = \delta(112, 122) = \frac{3}{4}$.
\end{cor}


\section{Average pattern co-occurrence}\label{avgSect}

In this section we deal with average, rather than maximal, pattern
co-occurrence.

Consider $S_n$ as a sample space with uniform distribution. Let
$\pi\in S_m$, and let $X_\pi$ be a random variable such that
$X_\pi(\tau)$ is the number of occurrences of pattern $\pi$ in a
given permutation $\tau\in S_n$.

It is an easy exercise to show that, even though the maximal number
of times a pattern can occur in a permutation (or a word, in
general) differs with the pattern, the average number of occurrences
of any pattern over all permutations of a given length is the same.

\begin{lem} \label{average}
$\displaystyle E(X_\pi)=\frac{1}{m!}\binom{n}{m}\sim
\frac{1}{(m!)^2}\,n^m$ for any pattern $\tau\in S_m$.
\end{lem}

\begin{proof}
Pick an $m$-letter subset $S$ of $[n]=\{1,2,\ldots,n\}$ in
$\binom{n}{m}$ ways. There is a unique permutation $\pi(S)$ of $S$
order-isomorphic to $\pi$, out of $m!$ equally likely permutations
in which the elements of $S$ can occur in $\tau\in S_n$. Let
$Y_{\pi(S)}$ be a random variable such that $Y_{\pi(S)}(\tau)$ is
the number of occurrences of $\pi(S)$ in $\tau$. Then
\[
P\!\left(Y_{\pi(S)}(\tau)=1\right)=\frac{1}{m!} \quad \text{and}
\quad P\!\left(Y_{\pi(S)}(\tau)=0\right)=1-\frac{1}{m!},
\]
so $E(Y_{\pi(S)})=1/m!$. This is true for any $S\subseteq[n]$ such
that $|S|=m$, and we have
$
X_\pi=\sum_{S\subseteq[n],\,|S|=m}{Y_{\pi(S)}},
$
hence,
\[
E(X_\pi)=\!\!\!\!\!\!\sum_{S\subseteq[n],\,|S|=m}{\!\!\!\!\!\!E(Y_{\pi(S)})}=
\frac{1}{m!}\binom{n}{m}.
\]
\end{proof}

Hence, the average pattern occurrence is the same. However, the
average pattern co-occurrence, measured by the covariance
$Cov(X_{\pi_1},X_{\pi_2})$, does depend on the pattern.
We will start by considering the average pattern co-occurrence with
itself, i.e.\ $V\!ar(X_\pi)$. That, via the standard deviation
$\sigma(X_\pi)$, will also tell us how tightly the distribution of
$X_\pi$ is grouped around the mean of Lemma \ref{average}.

Let $P_\pi$ be the permutation matrix of $\pi$, in other words,
$P_\pi=[\delta(\pi(i),j)]_{m\times m}$, where $\delta$ is the
Kronecker symbol. Note that $P_\pi$ is orthogonal, so
$P_{\pi^{-1}}=P_\pi^{-1}=P_\pi^t$. Also, for an integer $m>0$, and
integers $1\le i,j\le m$, define
\[
[i,j]_m=\binom{i-1+j-1}{i-1}\binom{m-i+m-j}{m-i}.
\]
Let $A_m$ be the $m\times m$ matrix with $(A_m)_{ij}=[i,j]_m$,
which have been studied e.g.\ in \cite{AE}.

\begin{thm} \label{variance} $\displaystyle
V\!ar(X_\pi)=c(\pi) n^{2m-1}+O(n^{2m-2})$ for any pattern $\pi\in
S_m$, $m>1$, where
\[
c(\pi)=\frac{1}{((2m-1)!)^2}
\left(Tr(A_mP_{\pi}A_mP_{\pi}^{-1})-\binom{2m-1}{m-1}^2\right)>0.
\]
The trace in the above formula can be expressed as
\[
Tr(A_mP_{\pi}A_mP_{\pi}^{-1})=\sum_{i,j=1}^{m}{[i,j]_m[\pi(i),\pi(j)]_m}
\]
For the standard deviation this gives
$\displaystyle \sigma(X_\pi)=\sqrt{c(\pi)}\
n^{m-\frac{1}{2}}+O(n^{m-1})$ for any pattern $\pi\in S_m$.
\end{thm}


\begin{proof}
Since $V\!ar(X_\pi)=E(X_\pi^2)-E(X_\pi)^2$, and the value of
$E(X_\pi)$ was determined in Lemma~\ref{average}, it remains only
to consider $E(X_\pi^2)$. We have
\[
E(X_\pi^2)=
E\left(\bigg(\sum_{S\subseteq[n],\,|S|=m}{\!\!\!\!\!\!Y_{\pi(S)}}\bigg)^2\right)=
\!\!\!\!\!\!\sum_{\substack{S_1,S_2\subseteq[n]\\|S_1|=|S_2|=m}}
{\!\!\!\!\!\!E\left(Y_{\pi(S_1)}Y_{\pi(S_2)}\right)}.
\]
Of course, $Y_{\pi(S_1)}Y_{\pi(S_2)}=1$ if and only if both
$\pi(S_1)$ and $\pi(S_2)$ are subsequences of $\tau$, otherwise,
$Y_{\pi(S_1)}Y_{\pi(S_2)}=0$.

Let $S=S_1\cup S_2$,  and $|S_1\cap S_2|=\ell$, so $|S|=2m-\ell$.
We can pick a subset $S\subseteq[n]$ in $\binom{n}{2m-\ell}$ ways.
Note that any such $S$ is order-isomorphic to
$[2m-\ell]=\{1,2,...,2m-\ell\}$. Hence, the number of permutations
$\rho(S)$ of $S$ such that $\rho\restriction_{S_1}\cong\pi$ and
$\rho\restriction_{S_2}\cong\pi$ is the same for any $S$ of cardinality
$2m-\ell$ and depends only on $m$ and $\ell$.

Therefore, $E(X_\pi^2)$ is a linear combination of
$\left\{\binom{n}{2m-\ell}\,\mid\,0\le\ell\le m\right\}$ with
coefficients that are polynomials in $m$ and $\ell$. The degrees
in $n$ of both $E(X_\pi^2)$ and $E(X_\pi)^2$ are $2m$, and the
coefficient of $n^{2m}$ in $E(X_\pi)^2$ is $1/(m!)^4$. On the
other hand, $S=S_1\cup S_2$, $|S|=2m$ and $|S_1|=|S_2|=m$ imply
that $S_1\cap S_2=\emptyset$, so $Y_{\pi(S_1)}$ and $Y_{\pi(S_2)}$
are independent, and hence
\[
P\left(Y_{\pi(S_1)}Y_{\pi(S_2)}=1\right)=P\left(Y_{\pi(S_1)}=1\right)
P\left(Y_{\pi(S_2)}=1\right)=\left(\frac{1}{m!}\right)^2.
\]
Let $[x^d]P(x)$ denote the coefficient of $x^d$ in a given
polynomial $P(x)$. Since there are $\binom{2m}{m}$ ways to partition
a set $S$ of size $2m$ into two subsets of size $m$, the coefficient
of $\binom{n}{2m}$ in $E(X_\pi^2)$ is $\binom{2m}{m}/(m!)^2$. Hence,
\[
[n^{2m}]E(X_\pi^2)=\frac{1}{(2m)!}\frac{1}{(m!)^2}\binom{2m}{m}=
\frac{1}{(m!)^4}.
\]
Thus $[n^{2m}]E(X_\pi^2)=[n^{2m}]E(X_\pi)^2$, so
$\deg_n(V\!ar(X_\tau))\le 2m-1$, and hence,
$[n^{2m-1}]V\!ar(X_\tau)\ge 0$.

We have
\[
\begin{split}
[n^{2m-1}]E(X_\pi)^2=[n^{2m-1}]\left(\frac{1}{m!}\binom{n}{m}\right)^2
&=\frac{2}{(m!)^2}\cdot[n^{m}]\binom{n}{m}\cdot[n^{m-1}]\binom{n}{m}=\\
&=\frac{2}{(m!)^2}\cdot\frac{1}{m!}\cdot\left(-\frac{\binom{m}{2}}{m!}\right)=
-\frac{m(m-1)}{(m!)^4}
\end{split}
\]
Similarly, the coefficient of $n^{2m-1}$ in the $\binom{n}{2m}$-term
of $E(X_\pi^2)$ is
\[
-\frac{\binom{2m}{2}}{(2m)!}\frac{1}{(m!)^2}\binom{2m}{m}=
-\frac{m(2m-1)}{(m!)^4},
\]
so we only need to find the coefficient of the
$\binom{n}{2m-1}$-term of $E(X_\pi^2)$.

As we noted before, all subsets $S\subseteq[n]$ of the same size
(in our case, of size $2m-1$) are equivalent, so we may assume
$S=[2m-1]=\{1,2,\ldots,2m-1\}$. We want to find the number of
permutations $\rho$ of $S$ such that there exist subsets
$S_1,S_2\subseteq S$ of size $m$ for which we have $|S_1\cap
S_2|=1$ (so $S_1\cup S_2=S$) and $\rho\restriction_{S_1}\cong\pi$
and $\rho\restriction_{S_2}\cong\pi$.

Suppose that we want to choose $S_1$ and $S_2$ as above, together
with their positions in $S$, in such a way that the intersection
element $e$ is in the $i$th position in $\pi(S_1)$ and the $j$th
position in $\pi(S_2)$ (of course, $1\le i,j\le m$). Then $e$
occupies position $(i-1)+(j-1)+1=i+j-1$ in $S$. Hence, there are
$\binom{i-1+j-1}{i-1}$ ways to choose the positions for elements
of $\pi(S_1)$ and $\pi(S_2)$ to the left of $e$, and
$\binom{m-i+m-j}{m-j}$ ways to choose the positions for elements
of $\pi(S_1)$ and $\pi(S_2)$ to the right of $e$. On the other
hand, both $\pi(S_1)$ and $\pi(S_2)$ are naturally
order-isomorphic to $\pi$, hence, under that isomorphism $e$ maps
to $\pi(i)$ as an element of $S_1$ and to $\pi(j)$ as an element
of $S_2$. Since $e$ is the \emph{unique} intersection element,
exactly $\pi(i)-1$ elements in $S_1$ and exactly $\pi(j)-1$
elements in $S_2$, all distinct, must be less than $e$, the rest
of the elements of $S$ must be greater than $e$, so we must have
$e=(\pi(i)-1)+(\pi(j)-1)+1=\pi(i)+\pi(j)-1$. There are
$\binom{\pi(i)-1+\pi(j)-1}{\pi(i)-1}$ ways to choose the elements
of $S_1$ and $S_2$ which are less than $e$, and
$\binom{m-\pi(i)+m-\pi(j)}{m-\pi(j)}$ ways to choose the elements
of $S_1$ and $S_2$ which are greater than $e$.

Thus, the positions of $\pi(e)$ in $\pi(S_1)$ and $\pi(S_2)$
uniquely determine the position $e$ and value $\pi(e)$ of the
intersection element; there are $[i,j]_m$ ways to choose which
other positions are occupied by $\pi(S_1)$ and which ones, by
$\pi(S_2)$; and there are $[\pi(i),\pi(j)]_m$ ways to choose which
other values are in $\pi(S_1)$ and which ones are in $\pi(S_2)$.

Now that we have chosen both positions and values of elements of
$S_1$ and $S_2$, we can produce a unique permutation $\rho(S)$ of
$S$ which satisfies our conditions above. Simply fill the positions
for $S_1$, resp. $S_2$, by elements of $\pi(S_1)$, resp.
$\pi(S_2)$, in the order in which they occur.

Since the total number of permutations of $S$ is $(2m-1)!$, the
coefficient of the $\binom{n}{2m-1}$-term of $E(X_\pi^2)$ is
\[
\frac{\sum_{i,j=1}^{m}{\binom{i-1+j-1}{i-1}\binom{m-i+m-j}{m-j}
\binom{\pi(i)-1+\pi(j)-1}{\pi(i)-1}\binom{m-\pi(i)+m-\pi(j)}{m-\pi(j)}}}
{(2m-1)!}
=\frac{\sum_{i,j=1}^{m}{[i,j]_{m}[\pi(i),\pi(j)]_{m}}}{(2m-1)!}.
\]
The coefficient of $n^{2m-1}$ in $V\!ar(X_\pi)$ is, by the previous
equations,
\[
\begin{split}
[n^{2m-1}]V\!ar(X_\pi)&=\frac{\sum_{i,j=1}^{m}{[i,j]_{m}[\pi(i),\pi(j)]_{m}}}{((2m-1)!)^2}-
\frac{m(2m-1)}{(m!)^4}+\frac{m(m-1)}{(m!)^4}\\
&=\frac{\sum_{i,j=1}^{m}{[i,j]_{m}[\pi(i),\pi(j)]_{m}}}{((2m-1)!)^2}-
\frac{1}{(m!(m-1)!)^2}\\
&=\frac{1}{((2m-1)!)^2}\left(\sum_{i,j=1}^{m}{[i,j]_{m}[\pi(i),\pi(j)]_{m}}
-\binom{2m-1}{m-1}^2\right).
\end{split}
\]
Since $c(\pi)$ is the leading coefficient of $V\!ar(X_\pi)$ (a
polynomial in $n$), we have $c(\pi)\ge 0$. The following lemma
implies that $c(\pi)>0$, which finishes the proof of Theorem
\ref{variance}.
\end{proof}

\begin{lem} \label{strict}
For any $\pi\in S_m$,
$\displaystyle\sum_{i,j=1}^{m}{[i,j]_{m}[\pi(i),\pi(j)]_{m}}>\binom{2m-1}{m-1}^2$.
\end{lem}

\begin{proof}
The matrix $A_m$ is symmetric and hence diagonalizable, and the
eigenvalues of $A$ are all distinct and known to be
$\{\binom{2m-1}{i-1}$, $1\le i\le m\}$. Each row of $A_m$
sums to $\binom{2m-1}{m-1}$, so $[1,1,\dots,1]$ is an eigenvector
for the largest eigenvalue $\binom{2m-1}{m-1}$. The same is true of
the similar matrix $P_{\pi}A_mP_{\pi}^{-1}$. Let
$D_m=[d_{ij}]_{m\times m}$ be the $m\times m$ diagonal matrix with
$d_{ii}=\binom{2m-1}{m-i}$ for $i=1,\dots,m$. Then $A_m=BD_mB^{-1}$
for some orthogonal matrix $B$, so recalling that $Tr(MN)=Tr(NM)$
for any $M,N$ for which $MN$ and $NM$ exist, we have
\[
\begin{split}
Tr(A_mP_{\pi}A_mP_{\pi}^{-1})&=Tr(BD_mB^{-1}P_{\pi}BD_mB^{-1}P_{\pi}^{-1})
=Tr(D_mB^{-1}P_{\pi}BD_mB^{-1}P_{\pi}^{-1}B)\\
&=Tr(D_m(B^{-1}P_{\pi}B)D_m(B^{-1}P_{\pi}B)^{-1})=Tr(D_mCD_mC^{-1}),
\end{split}
\]
and the matrix $C=B^{-1}P_{\pi}B=[c_{ij}]_{m\times m}$ is
orthogonal, i.e. $C^{-1}=C^t$. Let $\mathbf{b}_i$ be the $i$th
column of $B=[b_{ij}]_{m\times m}$. Then
$c_{ij}=\sum_{k=1}^{m}{b_{ki}b_{\pi(k)j}}$. In particular, the
column $\mathbf{b}_1=[1,\dots,1]^t$ remains unchanged for any
$\pi\in S_n$, so $c_{1j}=\mathbf{b}_1\cdot\mathbf{b}_j=\delta(1,j)$
and $c_{i1}=\mathbf{b}_i\cdot\mathbf{b}_1=\delta(i,1)$. Now we know
that
\[
Tr(D_mCD_mC^t)=\sum_{i,k=1}^{m}{d_{ii}d_{kk}c_{ik}^2},
\]
$c_{11}=1$, $c_{i1}=c_{1i}=0$ for $i>1$, so $c_{ij}$ are not all
zero for $i,j>1$ (otherwise $C$ is not invertible, let alone
orthogonal), hence
\[
Tr(D_mCD_mC^t)>\binom{2m-1}{m-1}^2.
\]
This proves the lemma.
\end{proof}

\begin{rem} \label{symmetries}
Note that the sum $\sum_{i,j=1}^{m}{[i,j]_{m}[\pi(i),\pi(j)]_{m}}$
is invariant under the symmetry operations on $S_m$: reversal
${\frak r}:i\mapsto m-i+1$, complement ${\frak c}:\pi(i)\mapsto
m-\pi(i)+1$, and inverse ${\frak i}:\pi\mapsto\pi^{-1}$.
Invariance under $\frak r$ and $\frak c$ also extends to
permutations of multisets. Thus permutations $\pi$ in the same
symmetry class $\bar\pi$ have the same $c(\pi)$. The
values of
\[
\Delta(\pi)=((2m-1)!)^2c(\pi)=
\sum_{i,j=1}^{m}{[i,j]_{m}[\pi(i),\pi(j)]_{m}}-\binom{2m-1}{m-1}^2
\]
for symmetry classes in $S_4$ ($m=4$) are given in the table below:
\begin{center}
\begin{tabular}{c|c|c|c|c|c}
  $\bar\pi$ & 1234 & 1243 & 1432 & 1342 & 2413 \\
  \hline
  $\Delta(\pi)$ & 491 & 359 & 327 & 239 & 91 \\
\end{tabular}
\end{center}
\end{rem}

\begin{rem} \label{least-delta}
It is easy to see that, for a given $m$, $\Delta(\pi)$ attains its
maximum when $\pi=id_m=12\dots m$ since the sequences
$\{[i,j]_{m}\}$ and $\{[\pi(i),\pi(j)]_{m}\}$ (with multiplicities)
are arranged in the same order. It would be interesting to
characterize the permutations $\pi_*$ for which
$\Delta(\pi_*)=\min_{\pi\in S_m}{\Delta(\pi)}$. For small values of
$m$, these permutations $\pi_*$ are:
\begin{center}
\begin{tabular}{c|c|c|c|c|c|c|c|c|c}
  $m$ & 1 & 2 & 3 & 4 & 5 & 6 & 7 & 8 & 9\\
  \hline
  $\pi_*(m)$ & 1 & 12 & 132 & 2413 & 25314 & 254163 & 3614725 & 37145826 & 385174926 \\
\end{tabular}
\end{center}
Interestingly, the patterns $\pi_*(m)$ are also less avoided than
most patterns of the same length, and in fact, are the least
avoided patterns for $m\le 5$.
\end{rem}

\bigskip

We can consider the co-occurrence of any two permutation patterns
similarly. Since the proof is similar to the variance case, it is
omitted.

\begin{thm} \label{covariance}
For any patterns $\pi_1,\pi_2\in S_m$, $m>1$, the covariance
$Cov(X_{\pi_1},X_{\pi_2})$ is given by
\[
Cov(X_{\pi_1},X_{\pi_2})=c(\pi_1,\pi_2) n^{2m-1}+O(n^{2m-2}),
\]
where
\[
c(\pi)=\frac{1}{((2m-1)!)^2}
\left(Tr(A_mP_{\pi_1}A_mP_{\pi_2}^{-1})-\binom{2m-1}{m-1}^2\right).
\]
The trace in the above formula is
\[
Tr(A_mP_{\pi_1}A_mP_{\pi_2}^{-1})=\sum_{i,j=1}^{m}{[i,j]_m[\pi_1(i),\pi_2(j)]_m}
\]
\end{thm}

Considering symmetry classes of pairs of patterns (see Remark
\ref{symmetries}), we see that there are 7 classes of pairs of
3-letter permutations: $\{123,123\}$, $\{132,132\}$,
$\{123,132\}$, $\{132,213\}$, $\{132,231\}$, $\{123,231\}$,
$\{123,321\}$ (listed in the order of decreasing asymptotical
covariance). The first two pairs obviously have a positive
covariance, but of the other five pairs, only $\{123,132\}$ has a
positive covariance.

It would be interesting to characterize the pairs
$\{\pi_1,\pi_2\}$ according to the sign or magnitude of their
covariance.

\bigskip

We now consider patterns contained in words, where repeated letters
are allowed both in the pattern and the ambient string.
The additional condition on a pattern $\pi\in[l]^m$ on
words, i.e.\ an a pattern of $m$ letters over an alphabet
$[l]=\{1,2,\dots,l\}$, is that $\pi$ must contain all letters in
$[l]$. We will also assume that the ambient strings are in the set
$[k]^n$.

\begin{thm} \label{variance-words}
Let $\pi$ be a map of $[m]=\{1,2,\ldots,m\}$ onto
$[l]=\{1,2,\ldots,l\}$. Then for any positive integers $l\le m$,
\[
V\!ar(X_\pi)=c(\pi)n^{2m-1}k^{2l-1}+O(n^{2m-2}k^{2l-1}+n^{2m-1}k^{2l-2}),
\]
where
\[
c(\pi)=\frac{1}{(2m-1)!(2l-1)!}\left(Tr(A_mP_{\pi}A_lP_{\pi}^{-1})
-\frac{(2m-1)!(2l-1)!}{((m-1)!)^2(l!)^2}\right).
\]
The trace in the above formula is
\[
Tr(A_mP_{\pi_1}A_lP_{\pi_2}^{-1})=\sum_{i,j=1}^{m}{[i,j]_m[\pi_1(i),\pi_2(j)]_l}.
\]
\end{thm}

\begin{rem} \label{rem-special}
Note also that, given $1\le l\le m$, Theorem \ref{variance-words}
applies to $l!S(m,l)$ patterns $\tau$, where $S(m,l)$ is the
Stirling number of the second kind.
\end{rem}

The proof of Theorem \ref{variance-words} is an obvious extension of
the proof of Theorem \ref{variance}. Unfortunately, the same
extension to words does not work for Lemma \ref{strict}, but only
yields
\[
\sum_{i,j=1}^{m}{[i,j]_m
[\pi(i),\pi(j)]_l}>\binom{2m-1}{m}\binom{2l-1}{l}
=\frac{l}{m}\left(\frac{(2m-1)!(2l-1)!}{((m-1)!)^2(l!)^2}\right),
\]
which is a weaker result than what we want.

There is a similar covariance result on words as well.

\begin{thm} \label{covariance-words}
For any patterns $\pi_1,\pi_2\in[l]^m$, $1<l\le m$, the covariance
$Cov(X_{\pi_1},X_{\pi_2})$ is
\[
Cov(X_{\pi_1},X_{\pi_2})=c(\pi_1,\pi_2)n^{2m-1}k^{2l-1}
+O(n^{2m-2}k^{2l-1}+n^{2m-1}k^{2l-2}),
\]
where
\[
c(\pi_1,\pi_2)=\frac{1}{(2m-1)!(2l-1)!}\left(Tr(A_mP_{\pi_1}A_lP_{\pi_2}^{-1})
-\frac{(2m-1)!(2l-1)!}{((m-1)!)^2(l!)^2}\right).
\]
The trace in the above formula is
\[
Tr(A_mP_{\pi_1}A_mP_{\pi_2}^{-1})=\sum_{i,j=1}^{m}{[i,j]_m[\pi_1(i),\pi_2(j)]_l}.
\]
\end{thm}


\end{document}